\documentclass[final,nomarks]{dmtcs-episciences}

\usepackage[utf8]{inputenc}
\usepackage[english]{babel}
\usepackage{amsmath,amssymb,amsfonts,mathtools}
\usepackage[numbers]{natbib}
\hypersetup{hidelinks}

\newcommand{\bb}{\begin{equation}}
\newcommand{\ee}{\end{equation}}

\newcommand{\QED}{\hfill$\square$}

\newtheorem{thm}{Theorem}
\newtheorem{conj}[thm]{Conjecture}
\newtheorem{prop}[thm]{Proposition}

\title{The transversal game}
\author{\v{Z}arko Ran\dj elovi\'c}
\affiliation{Mathematical Institute of the Serbian Academy of Sciences and Arts, Belgrade, Serbia}

\begin{document}
\publicationdata{vol. 28:3}{2026}{26}{10.46298/dmtcs.19063}{2026-08-27; 2026-08-27; 2026-11-17}{2026-11-19}
\maketitle

\begin{abstract}
    Two players take turns claiming empty cells from an $n\times n$ grid. The first player (if any) to occupy a transversal (a set of $ n $ cells having no two cells in the same row or column) is the winner. What is the outcome of the game given optimal play? Our aim in this paper is to show that for $n\ge 4$ the first player has a winning strategy. This answers a question of Erickson. We also consider variants of the game in which the family of winning transversals
    is restricted, both in the original game and in its
    Maker--Breaker version.
\end{abstract}
\section{Introduction}

Suppose that player $1$ and player $2$ play a game on an $n\times n$ grid. They play alternating moves. Player $1$ goes first. Each player when it is their move writes their letter (say $X$ for player $1$ and $O$ for player $2$) in one of the unoccupied cells of the grid. The first player to occupy a set of $n$ cells having no two cells in the same row or column (meaning that those $n$ cells only have the letter of that particular player) is the winner. If neither player does that and the grid is filled the game is a draw. Such positional games have been studied widely (some examples include \cite{beck}, \cite{erdos}, \cite{hales}). Our aim is to prove that player $1$ wins if $n\ge 4$ (it
is easy to check that it is a draw for $n=2$ and $n=3$; see Proposition \ref{3draw}). This question was posed by Erickson \cite{transversal}. We call the above described game the \textit{transversal $n$-game} where $n\times n$ is the size of the grid.

The transversal $n$-game is a strong positional game in the usual terminology (see, for example, \cite{positional}). Its board is the set of $n^2$ cells, and its winning sets are the transversals. If the rows and columns are identified with the two vertex classes of the complete bipartite graph $K_{n,n}$, each cell $(i,j)$ corresponds to an edge and each transversal corresponds to a perfect matching. Thus our game is precisely the strong perfect-matching game on $K_{n,n}$. By the standard strategy-stealing argument for strong positional games (see, for example, \cite{hales,positional}), player $1$ can always guarantee at least a draw from the empty board. Consequently, the question is whether player $1$ has a winning strategy or player $2$ can force a draw.

\begin{thm}\label{maintransversal}
    For $n=1$ and $n\ge 4$ player $1$ wins the transversal $n$-game.
\end{thm}

We say that a \textit{transversal} is a set of $n$ cells on an $n\times n$ grid having no two cells in the same row or column. We will also use $[k]$ to denote the set $\{1,2,..,k\}$ for any $k$.

When considering an $n\times n$ grid we will label the cells of the grid with ordered pairs $(a,b)$ with $1\le a,b\le n$. The cell $(a,b)$ will be the cell in the $a$-th row from top to bottom and the $b$-th column from left to right (similar to the labeling in matrices).

We give an example of a game for $n=4$. Player $X$ (i.e. player $1$) makes his first move in the cell $(1,1)$. Then player O (i.e. player $2$) makes his first move in $(2,2)$. Then player $X$ plays in $(2,3)$ trying to make a threat on the next move. Then player $O$ plays in $(4,4)$ which blocks a potential threat and gives $O$ a possible chance to make his own threat on the next move depending on what $X$ plays. Now $X$ plays in $(4,2)$ making a threat on $(3,4)$. Now, $O$ defends by playing in $(3,4)$. Now $X$ plays in $(3,1)$ making a threat on $(1,4)$. Since $O$ has no threat, $O$ must play in $(1,4)$. Now $X$ plays in $(2,4)$ and makes a threat on both $(1,3)$ and $(3,3)$ (labeled $W$ below). Since $O$ still has no threats and cannot defend both of $X$'s threats, $X$ wins on the next move. Subscripts below indicate which move it is for that player.
$$\begin{pmatrix}
X_1 &   & W & O_4\\
 & O_1 & X_2  & X_5\\
X_4 &   & W  & O_3 \\
 & X_3 &  & O_2
\end{pmatrix}$$

We give a brief outline of the strategy. In the proof of Theorem 1 we will have player $X$ making a threat as early as possible. From then on the key point is for player $X$ to make moves that help to achieve his goal but without allowing player $O$ to make threats. Player $X$ will start making threats as soon as he has occupied $n-1$ cells  of the grid and will continue to make new threats forcing $O$ to block those threats.

In Section 3 we also consider variants of the game in which
only transversals belonging to a prescribed family are
winning sets. We obtain bounds on the number of transversals
needed to guarantee a win for player $1$, and we also
investigate the corresponding Maker--Breaker version.

In the weaker Maker-Breaker version of this game (where $X$ wins if he occupies the cells of a transversal and $O$ wins as long as $X$ does not make a transversal) $X$ wins if $n\ge 4$. Indeed, this may be proved by induction. If $X$ can win for $n\ge 4$ while making the first move in the top-left cell then if we have an $(n+1)\times (n+1)$ grid suppose $X$ plays on $(1,1)$. If $O$ plays in row $1$ or column $1$ then $X$ can win by induction after removing the first row and column. Otherwise if $O$ plays in $(a,b)$ with $a,b>1$ then we can have $X$ playing in $(a,c)$ for $c\neq 1,b$. Now after removing row $a$ and column $c$ $X$ can win by induction since he has already occupied the top-left cell. One can check that $X$ can win with starting move $(1,1)$ for $n=4$ by a case analysis.

In a strong positional game both players try to complete a winning set before their opponent, while in a Maker--Breaker (or weak) game Maker tries to complete such a set and Breaker only tries to prevent this. A winning strategy for Maker does not by itself give a winning strategy for the first player in the corresponding strong game, since the opponent might complete a winning set first. This is the reason why controlling $O$'s threats is an essential part of our strategy.

Perfect-matching games have also been studied extensively on other graphs. Chv\'{a}tal and Erd\H{o}s \cite{chvatalerdos} initiated the systematic study of biased positional games on complete graphs. Fast strategies for the Maker--Breaker perfect-matching game on complete graphs were studied by Hefetz, Krivelevich, Stojakovi\'{c} and Szab\'{o} \cite{fastpm}, and biased versions were further investigated by Mikala\v{c}ki and Stojakovi\'{c} \cite{biasedfastpm}. Ben-Shimon, Ferber, Hefetz and Krivelevich \cite{randompm} studied perfect-matching Maker--Breaker games on random graphs, obtaining a hitting-time result for the random graph process. These papers concern Maker--Breaker games and do not directly establish a strong-game win on $K_{n,n}$. Ferber and Hefetz \cite{ferberhefetz} obtained explicit winning strategies for the strong perfect-matching game on complete graphs of sufficiently large even order by drawing on fast winning strategies for the corresponding weak game. Their complete-graph setting is different from our bipartite one: a strategy on $K_{2n}$ need not use only the edges of $K_{n,n}$. Our proof instead gives a forcing strategy directly on the bipartite board, while ensuring that $O$ cannot make a winning threat at the critical stages.

There were two earlier arXiv submissions concerning Erickson's question: Krishna \cite{krishna} withdrew his proposed proof after identifying an error, and Dumitrescu \cite{dumitrescu} withdrew his analysis because the rules had been misinterpreted. Neither withdrawn submission constitutes an established solution. After the first version of the present paper appeared, we learned of an independent proof by Guan \cite{guan}, which also gives a strategy forcing a win by player $1$'s $(n+2)$-nd move.

\section{Proof of Theorem 1}
 We will say that the \textit{value} of a position is `player $1$' if the game is a player $1$ win, `player $2$' if a player $2$ win, and `draw' if it is a draw (all with
perfect play of course). Also we say that the \textit{value} of a game is the value of its starting position. Permuting rows or columns also permutes the transversals. This means that permuting rows or columns does not change the value of a position. Notice that reflection with respect to the main diagonal of the grid (i.e. cells $(i,i)$) also does not change the value of a position. We define \textit{good} transformations of the grid to be any combination of reflecting with respect to the main diagonal and permuting rows or columns.

 Let us say that regardless of whose move it is a player has a \textit{threat} on cell $(a,b)$ if $(a,b)$ is empty and they would win on their next turn (if they have one) by playing on $(a,b)$. If a player can not win on their next move we say that player has no threats.

For completeness, and to give the reader a feel for the game, we show
that for $n=3$ the value is a draw. 
\begin{prop}\label{3draw}
 The value of the transversal 3-game is a draw.  
\end{prop}

\textit{Proof.} We will show that both $X$ and $O$ can guarantee at least a draw. It is easy to see that $X$ can guarantee at least a draw by a strategy stealing argument, similar to Lemma 4 in \cite{hales}.

Now we show that $O$ can guarantee a draw. By permuting rows and columns every opening move for $X$ is equivalent, so we may assume that $X$ plays on $(1,1)$. Then $O$ plays on $(2,2)$. Now by symmetry with respect to the main diagonal we may assume that $X$ plays on either $(1,2),(1,3),(2,3)$ or $(3,3)$. We deal with each case individually.

If $X$ plays on $(1,2)$ then $O$ plays on $(2,3)$. Now if $X$ does not play on $(2,1)$ then $O$ plays there next which blocks off a row, but if $X$ does play on $(2,1)$ then $O$ plays on $(3,3)$. Now again $X$ has no threats and if $X$ does not play on $(1,3)$ then $O$ can block off a column, but if $X$ does play on $(1,3)$ then $O$ plays on $(3,2)$ which blocks all transversals. 

$$\begin{pmatrix}
X_1 & X_2 & X_4\\
X_3 & O_1 & O_2\\
 & O_4 & O_3
\end{pmatrix}$$

If $X$ plays on $(1,3)$ then $O$ plays on $(2,3)$. Now if $X$ does not play on $(2,1)$ then $O$ plays there next which blocks off a row, but if $X$ does play on $(2,1)$ then $O$ plays on $(3,2)$. Now again if $X$ does not play on $(1,2)$, then $O$ can block off a column, but if $X$ does play on $(1,2)$ then $O$ plays on $(3,3)$ which blocks all transversals.

$$\begin{pmatrix}
X_1 & X_4 & X_2\\
X_3 & O_1 & O_2\\
 & O_3 & O_4
\end{pmatrix}$$

If $X$ plays on $(2,3)$ then $O$ plays on $(3,2)$. Now if $X$ does not play on $(1,2)$ then $O$ plays there next which blocks off a column, but if $X$ does play on $(1,2)$ then $O$ plays on $(3,1)$. Now again if $X$ does not play on $(3,3)$ then $O$ can block off a row, but if $X$ does play on $(3,3)$ then $O$ plays on $(2,1)$ which blocks all transversals.

$$\begin{pmatrix}
X_1 & X_3 & \\
O_4 & O_1 & X_2\\
O_3 & O_2 & X_4
\end{pmatrix}$$

Finally if $X$ plays on $(3,3)$ then $O$ plays on $(3,2)$. Now if $X$ does not play on $(1,2)$ then $O$ plays there next which blocks off a column, but if $X$ does play on $(1,2)$ then $O$ plays on $(2,1)$. Now again if $X$ does not play on $(2,3)$ then $O$ can block off a row, but if $X$ does play on $(2,3)$ then $O$ plays on $(3,1)$ which blocks all transversals.

$$\begin{pmatrix}
X_1 & X_3 & \\
O_3 & O_1 & X_4\\
O_4 & O_2 & X_2
\end{pmatrix}$$

So player $O$ can always guarantee at least a draw.\hfill\QED

Now we move onto Theorem $1$ where $n\ge 4$.

\textit{Proof of Theorem 1.} The case $n=1$ is trivial so suppose that $n\ge 4$. The general strategy of the proof will be as follows. After the first $2n-2$ moves in total (up to a permutation of rows and columns) $X$ will occupy the cells of the main diagonal except for the bottom right cell and $O$ will occupy at least two cells from the last column, one of which will be forced to be $(n,n)$. Then $X$ will make threats and will only need up to four more moves to win depending on where $O$ has played.

Player $X$ starts by playing on $(1,1)$. Now we will show by induction that player $X$ can guarantee that the following holds after $k+1$ $X$-moves and $k$ $O$-moves for all $1\le k \le n-3$: Up to some good transformations we have that:

\begin{itemize}

 \item There are $X$s in cells $(a,a)$ for $1\le a\le k+1$,\\[-4ex]
 \item There are no $O$s in cells $(a,b)$ with $a,b>k+1$,\hfill (*)\\[-4ex]
 \item There is at least one cell $(a,b)$ with an $O$ with $b=k+2$  
\end{itemize}

We first show this for $k=1$. Since $X$ played in $(1,1)$ we may assume by applying good transformations that $O$ played either in $(2,2)$ or $(1,2)$. In either case $X$ will play in $(2,3)$ and we will swap columns $2$ and $3$. Now we have that $(1,1),(2,2)$ have $X$s and either $(1,3)$ or $(2,3)$ has an $O$. This shows the case $k=1$.

Now suppose that we can do this for $1\le k\le n-4$. So after $k+1$ $X$-moves and $k$ $O$-moves (*) is satisfied. Suppose now that $O$ plays on cell $(a,b)$.

If either $a\le k+1$ or $b\le k+1$ then $X$ will play in $(k+2,k+3)$ and we will swap columns $k+2$ and $k+3$. This will ensure (*) is satisfied after $k+2$ $X$-moves and $k+1$ $O$-moves.

Suppose now that $a,b>k+1$. Now if $b>k+2$ then $X$ will play in $(a,k+2)$ and we will swap rows $a$ and $k+2$ and then if $b\neq k+3$ swap columns $k+3$ and $b$. If $b=k+2$ then we will swap columns $k+2$ and $k+3$. Now we do the same as in the case with $b>k+2$. In either case we ensure that (*) is satisfied which completes the induction step.

$$\begin{pmatrix}
X & & & & &\\
 &\ddots & &O & &\\
 & &X & & & \\
 & & & & &\\
 & & & & &\\
 & & & & &
 \end{pmatrix}$$

By above player $X$ can guarantee that (*) is satisfied for $k=n-3$. After $n-2$ $X$-moves and $n-3$ $O$-moves cells $(n-1,n),(n,n)$ are empty. Now if $O$ plays on $(n,n-1)$ or $(n-1,n)$ $X$ will play on $(n,n)$, otherwise $X$ will play on $(n-1,n)$. Either way $X$ will have a threat on $(n-1,n-1)$ or $(n,n-1)$. Since $O$ has no threats as there are only $n-2$ $O$s so far $O$ must defend by blocking $X$'s threat. After a total of $2n-2$ moves by swapping columns $n-1$ and $n$ and if needed swapping rows $n-1$ and $n$ by (*) we can make sure that $X$ occupies cells $(i,i)$ for $1\le i\le n-1$ while $O$ occupies cell $(n,n)$, another cell $(a,n)$ for $a<n$ and some other cells. At this point after a total of $2n-2$ moves let $T_O,T_X$ be the sets of cells with an $O,X$ respectively. Now let $$S_O=\{i\in [n]|(i,n)\not \in T_O \wedge (n,i)\not \in T_O\}$$
Since $|T_O|=n-1$ we know that $|S_O|\ge 1$. But also since $(n,n),(a,n)\in T_O$ we know that $|S_O|\le n-2$. Now notice that since there are at least two $O$s in the final column, there are at least 2 columns without $O$s. This means that $O$ currently has no threat. By using the same permutation on the first $n-1$ rows and the first $n-1$ columns we can guarantee that $S_O=[s]$ for some $s\le n-2$. For a set $A\subset [n]$ we denote $T_A=\{(i,i)|i\in[n]\backslash A\}$. Now we split into three cases depending on $s$ and $n$.

Case 1: $s=1,n\ge 5$. Let $r=|\{i<n|(i,n)\in T_O\}|$. By above $1\le r\le n-2$. Since each pair of cells $\{(i,n),(n,i)\}$ for $2\le i\le n-1$ has an $O$ in at least one cell and $(n,n)$ has an $O$ but $|T_O|=n-1$, we must have that every such pair has exactly one $O$. In fact these together with $(n,n)$ must be all the $O$s in the grid. Now similar to above by using the same permutation on rows and columns we may assume that without loss of generality $$T_O=\{(n,i)|2\le i\le n-r-1\}\cup \{(i,n)|n-r\le i\le n\}.$$ Now $X$ plays on $(1,n)$ and since $(1,n),(n,1),T_{\{1,n\}}$ form a transversal, $X$ makes a threat on $(n,1)$. So $O$ must defend the threat on $(n,1)$. Since all $O$s are in row $n$ and column $n$, $O$ has no threats. Now $X$ plays on $(n,n-1)$ and since $(n,n-1),(1,n),(n-1,1),T_{\{1,n-1,n\}}$ form a transversal $O$ must defend on $(n-1,1)$. Since all $O$s are in columns $1,n$ and row $n$ $O$ has no threats. This is because to have a threat $O$ must have at least $n-1\ge 4$ cells no two of which share a row or column. Now we split this into two subcases.

Case 1a: $r=n-2$. Now $X$ plays on $(2,1)$. We have that cells $(2,1),(1,n),(n,2)$ and $T_{\{1,2,n\}}$ form a transversal, but also $(n-1,2),(2,1),(1,n),(n,n-1)$ together with $T_{\{1,2,n-1,n\}}$ form a transversal. Thus $X$ has a threat on both $(n-1,2)$ and $(n,2)$. Since $O$ cannot defend both threats, $X$ will win on the next move.

$$\begin{pmatrix}
X &   &   &   &   & X_n\\
X_{n+2}& X &   &   &   & O\\
 & & X &   &   & O\\
  &   &   &\ddots&   &\vdots \\
O_{n+1}  & W  &   &   & X & O\\
O_n&W & &   & X_{n+1}  & O_{n-1}
\end{pmatrix}$$

Case 1b: $1\le r\le n-3$. Now $X$ plays on $(n-1,2)$. Since $(n,n-1),(n-1,2),(2,n),T_{\{2,n-1,n\}}$ form a transversal and $(2,1),(1,n),(n,n-1),(n-1,2),T_{\{1,2,n-1,n\}}$ form a transversal we have that $X$ has threats on both $(2,1)$ and $(2,n)$. Since $O$ has no threats, $X$ will win on the next move.

$$\begin{pmatrix}
 X& & & & & &  &X_n\\
W &X & & & & &  &W\\
 & &\ddots & & & & &\\
 & & & X& & &  &\\
 & & & &X & & & O\\
 & & & & &\ddots &  &\vdots\\
 O_{n+1}&X_{n+2} & &  & & &X &O\\
 O_n& O&\dots & O& &  &X_{n+1} &O_{n-1}\\
\end{pmatrix}$$

Case 2: $2\le s\le n-2$. Since there is at least one $O$ in the cells $(i,n),(n,i)$ for all $s+1\le i\le n$ the number of $O$s in the first $s$ columns is at most $(n-1)-(n-s)=s-1$. Since $s\ge 2$ the number of ordered pairs $(i,j)$ such that $1\le i,j\le s$ and $i\neq j$ is $s(s-1)>s-1$. Thus there are $b,c\in [s]$ such that $b\neq c$ and $(b,c)$ is an unoccupied cell. Now $X$ plays in $(n,b)$. Since $(n,b),(b,n),T_{\{b,n\}}$ form a transversal $O$ must play in $(b,n)$. Since there were already $O$s in the final column and there were at least $2$ columns without $O$s we still have $2$ columns without $O$s and hence $O$ has no threats. Now $X$ plays on $(c,n)$. Since $(c,n),(n,c),T_{\{c,n\}}$ form a transversal and $(n,b),(b,c),(c,n),T_{\{b,c,n\}}$ form a transversal, $X$ has threats on both $(n,c)$ and $(b,c)$. This means that $X$ will win on the next move.

$$\begin{pmatrix}
 X& & & & & & & \\
 &\ddots & & & & & & \\
 & &X & &W & & &O_n \\
 & & &\ddots & & & & \\
 & & & &X & & &X_{n+1} \\
 & & & & &\ddots & & \\
 & & & & & &X & \\
 & &X_n & &W & & &O_{n-1} 
\end{pmatrix}$$

Case 3: $s=1,n=4$. Since column $4$ has at least two $O$s by using the same permutation of rows and columns we can assume that there is an $O$ in $(3,4)$. Now $X$ plays on $(1,4)$ which forces $O$ to defend on $(4,1)$. Now $O$ has no threats since all $O$s are in row $4$ and column $4$ and to have a threat $O$ would need $3$ cells no two of which share a row or column. Since $s=1$ and $(3,4)$ has an $O$ we know that exactly one of $(2,4)$ and $(4,2)$ have an $O$. Thus $X$ can play on $(4,3)$. Since $(3,1),(1,4),(4,3),(2,2)$ is a transversal, $X$ has a threat on $(3,1)$ and hence $O$ must defend on $(3,1)$. We split into two subcases:

Case 3a: If $(2,4)$ has an $O$.  Notice that $O$ still has no threat and $X$ will play on $(2,1)$. Now since $(2,1),(1,4),(4,3),(3,2)$ is a transversal and also $(2,1),(1,4),(4,2),(3,3)$ is a transversal, $X$ has threats on both $(3,2)$ and $(4,2)$ and will hence win on the next move.

$$\begin{pmatrix}
X &   &  & X_4\\
X_6& X &   & O\\
O_5& W  & X  & O \\
 O_4&W & X_5 & O_3
\end{pmatrix}$$

Case 3b: If $(4,2)$ has an $O$ we still have that $O$ has no threats so $X$ can play on $(3,2)$. Since $(3,2),(2,4),(4,3),(1,1)$ is a transversal and $(3,2),(2,1),(1,4),(4,3)$ is also a transversal, $X$ has threats on both $(2,1)$ and $(2,4)$ and will hence win on the next move.

$$\begin{pmatrix}
X &   &  & X_4\\
W & X &   & W \\
O_5& X_6  & X  & O \\
 O_4& O & X_5 & O_3
\end{pmatrix}$$

This covers all cases and thus we have shown that for all $n\ge 4$ player $X$ has a winning strategy. This concludes the proof. \hfill\qed
\section{Related questions}
We now know that player $X$ wins for $n\ge 4$ so a natural question to ask is how many transversals are needed for player $X$ to win.
Suppose that we have a family $\mathcal{F}$ of transversals on an $n\times n$ grid and suppose that the goal of the game for both players was to occupy $n$ cells of some transversal in $\mathcal{F}$. We will call this game the transversal $\mathcal{F}$-game.

It is unclear if our strategy is optimal in terms of how many transversals are needed to guarantee that player $X$ wins. So one might ask what is the minimal possible $f(n)$ such that there is a family $\mathcal{F}$ with $|\mathcal{F}|=f(n)$ where player $X$ has a winning strategy in the transversal $\mathcal{F}$-game?

Another natural question to ask is what is the minimal number $g(n)$ such that for every family $\mathcal{F}$ of at least $g(n)$ transversals player $X$ wins?

It is clear that $f(n)\le g(n)\le n!$. We will give an argument why $f(n)<\frac{p(n)n!}{2^n}$ for some polynomial $p$ and also why $g(n)<n!-(n!)^{1/5-O(1/n)}$. We will be using a very similar strategy as in the proof of Theorem \ref{maintransversal} but being a bit more careful when choosing the first $n-1$ moves. 
\begin{thm}\label{thm3}
    There is a fixed polynomial $p$ such that $f(n)<\frac{p(n)n!}{2^n}$ for every $n\ge 4$.
\end{thm}
\begin{proof}
    Let $k=\lfloor \frac{n-2}{2}\rfloor$. We may assume that $n\ge 8$ since we can otherwise just increase $p$ to make our statement true for all $n\ge 4$. We choose arbitrary sets of $k$ rows and $k$ columns with the only restriction that one of the chosen rows is the 1st row and none of the chosen columns is the 1st column. Let $A,B\subset [n]$ be the chosen sets of rows and columns respectively. Consider the family $\mathcal{F}$ of all transversals that have at most four cells in $A\times B$. Our goal is to show that player $X$ can win the transversal $\mathcal{F}$-game. We will consider $A,B$ as variable sets so whenever we permute rows or columns or potentially perform a reflection with respect to the main diagonal we also accordingly change $A,B$ (and of course $\mathcal{F}$). For convenience denote by $A^i,B^i$ the sets $A,B$ after exactly $2i$ total moves. For clarity we will only perform good transformations after $X$'s move (before $O$ makes his next move) just like in the proof of Theorem \ref{maintransversal}. Since both $A,B$ have size $k$ the reflection with respect to the main diagonal will just swap those sets and not change their sizes. We first prove the following claim.

    \textbf{Claim.} If $X$ can achieve the following after the first $2n-2$ moves in the transversal $\mathcal{F}$-game:
    
    \begin{itemize}
        \item  $X$ occupies precisely the cells $(i,i)$ for $1\le i\le n-1$, 
        \item $O$ occupies $(n,n)$ and at least one other cell in the $n$-th column,
        \item none of the cells $X$ occupies are in $A\times B$,
    \end{itemize}then $X$ can guarantee a win.

    \textit{Proof of claim.} Suppose that $X$ achieves the required conditions after $2n-2$ total moves. Now $X$ follows the exact same strategy as in the proof of Theorem \ref{maintransversal}. Notice that in each case $X$ requires at most four extra moves to win. Since none of the first $n-1$ cells taken by $X$ are inside $A\times B$ we have that all the transversals that $X$ makes threats on prior to winning have at most four cells in $A\times B$ and hence belong to $\mathcal{F}$. Thus $X$ can win the transversal $\mathcal{F}$-game which proves the claim. \QED

    Now we will show that $X$ can indeed achieve the conditions required by the claim. After that we will show that $|\mathcal{F}|\le \frac{p(n)n!}{2^n}$ for some fixed polynomial $p$.

    For $1\le i\le 2k$ let $A_i=\{i+1,\ldots ,n\}\cap A^i$ and $B_i=\{i+1,\ldots ,n\}\cap B^i$. Notice that after a total of $2i$ moves $X$ needs to avoid $A_i\times B_i$ on his next move. 
    Now we start the game. Player $X$ begins by playing on $(1,1)$. Then $O$ plays somewhere. At this point since $1\in A\setminus B$ we have that $|A_1|=k-1, |B_1|=k$. We will show by induction that player $X$ can guarantee that the following holds after $i$ $X$-moves and $i$ $O$-moves for all $2\le i \le 2k$: 

\begin{itemize}

\item There are $X$s in cells $(x,x)$ for $1\le x\le i$,\\[-4ex]
 \item None of the first $i-1$ moves by $O$ are in cells $(x,y)$ with $x,y>i$,\hfill (**)\\[-4ex]
 \item There is at least one cell $(x,y)$ with an $O$ with $y=i+1$  
 \item $|A_i|+|B_i|\le 2k-i$ or $|A_i||B_i|=0$
\end{itemize}
We first show this for $i=2$. Suppose that $O$ plays his first move in some $(a,b)$. We first suppose that $b\neq 1$. There are fundamentally two subcases to consider depending on whether $a=1$ or $a\neq 1$. We will show that in either case $X$ can choose some $c\neq 1,b$ and play in $(\max(a,2),c)$ where we have either $\max(a,2)\in A_1$ or $c\in B_1$ but not both. Indeed if $\max(a,2)\in A_1$ then since $k<n-2$ pick $c\neq 1,b$ such that $c\not \in B_1$. If $\max(a,2)\not \in A_1$ then pick $c\in B_1$ with $c\neq 1,b$ which is possible since $k>2$. Now we swap columns $2$ and $c$ and if $a\neq 1$ swap rows $2$ and $a$. If $b=2$ then we further swap columns $c$ and $3$ and if $b\neq 2$ then we swap columns $b$ and $3$. Now we have that $(1,1),(2,2)$ have $X$s and either $(1,3)$ or $(2,3)$ has an $O$. 
After $O$'s next move we also have that $|A_2|+|B_2|\le |A_1|+|B_1|-1\le 2k-2$ ensuring that (**) is satisfied in the case $i=2$.

Now suppose that $b=1$. Since $(1,1)$ is already occupied by $X$, we have that $a\neq 1$. Choose $x\notin A_1\cup\{1,a\}$ and $c\in B_1$, and let $X$ play in $(x,c)$. Such a choice is possible since $k<n-2$ and $|B_1|>0$. Now by permuting rows and columns, fixing the first row and
column, and sending rows $x,a$ to rows $2,3$, respectively, and column $c$ to column $2$, we can ensure that $X$ occupies $(1,1)$ and $(2,2)$ while $O$ occupies $(3,1)$. Reflecting with respect to the main diagonal, we obtain an $O$ in $(1,3)$. Note that under this reflection the sets $A$ and $B$ are also interchanged. Since $x\notin A_1$ and $c\in B_1$, after $O$'s next move we have that $|A_2|+|B_2|=2k-2$. Thus $(**)$ is also satisfied in this case.

Now suppose that we can do this for $2\le i\le 2k-1$. So after $i$ $X$-moves and $i$ $O$-moves (**) is satisfied. Suppose now that $O$'s last move was in cell $(a,b)$.

If $|A_i||B_i|=0$ then we do the same step as in Theorem \ref{maintransversal} and we have that $(**)$ is satisfied.

Now suppose that $|A_i||B_i|\neq 0$ and either $a\le i$ or $b\le i$. Choose any $x\in A_i$. We need to choose $c\ge i+2$ such that $c\notin B_i$. Since $|B_i|\le 2k-i$ and $n-i-1>2k-i$ there is a $c\not \in B_i$ with $c\ge i+2$. Now $X$ will play in $(x,c)$. Then we will swap rows $x$ and $i+1$, columns $i+1$ and $c$ and if $c\neq i+2$ we will swap columns $i+2$ and $c$. This will ensure (**) is satisfied after $i+1$ $X$-moves and $i+1$ $O$-moves.

Suppose now that $|A_i||B_i|\neq 0$ and $a,b>i$. First suppose that $a\in A_i$. Now since $n-i-1>2k-i\ge |B_i|$ there is a $c\ge i+1,c\neq b$ such that $c\not \in B_i$. Now $X$ plays in $(a,c)$. Swap rows $a$ and $i+1$ and swap columns $c$ and $i+1$. Next, if $b=i+1$ swap columns $c$ and $i+2$, but if $b>i+1$ swap columns $b$ and $i+2$. Either way, the conditions of (**) are going to be satisfied after $2i+2$ total moves. Now suppose that $a\not \in A_i$. If there is a $c\in B_i\setminus \{b\}$ then we will have $X$ play in $(a,c)$. Now we perform the same swaps as in the above case (when $a\in A_i$) and achieve (**) after $2i+2$ total moves. Since $|B_i|\neq 0$ the only other case is if $B_i=\{b\}$. Now since $n-i-1>2k-i\ge |A_i|$ there is an $a'\ge i+1,a'\not \in A_i\cup \{a\}$. Notice that $a'\neq a$ so we will have $X$ play in $(a',b)$. Now swap columns $b$ and $i+1$ and swap rows $a'$ and $i+1$. Further if $a=i+1$ swap rows $a'$ and $i+2$, but if $a\neq i+1$ swap rows $a$ and $i+2$. Now perform a reflection with respect to the main diagonal. This ensures that (**) is satisfied after $2i+2$ total moves which completes the induction step.

After $2k$ $X$ moves and $2k$ $O$ moves we obtain that $A_{2k}=\emptyset$ or $B_{2k}=\emptyset$ and hence we can continue as in Theorem \ref{maintransversal} to occupy all cells $(i,i)$ with $1\le i\le n-1$. This ensures that $X$ has avoided $A\times B$ throughout all of his first $n-1$ moves. Note that this requires at most two more moves since $n-1-2k\le 2$. On $O$'s $(n-1)$-th move he must play in $(n,n)$ since at that point the main diagonal is a transversal in $\mathcal{F}$. This will ensure the conditions of the above claim meaning that $X$ can win the transversal $\mathcal{F}$-game.

Now we will estimate $|\mathcal{F}|$. For any $0\leq m\leq 4$,
let $P_{n,m}$ be the number of transversals that have
exactly $m$ cells in $A\times B$.
First we choose the $m$ rows from $A$ in which the transversal has a cell in $A\times B$. Then, from top to bottom, we choose $m$ distinct columns from $B$ for these cells. Next, again from top to bottom, we choose distinct columns outside $B$ for the remaining $k-m$ rows of $A$. Finally, we choose the cells in the remaining $n-k$ rows and $n-k$ columns, which can be done in $(n-k)!$ ways. Thus
 
    $$P_{n,m}=\binom{k}{m}k(k-1)\cdots (k-m+1)(n-k)(n-k-1)\cdots (n-2k+m+1)\cdot(n-k)!.$$ We have that 
    \begin{align}\label{pnkm}
    P_{n,m}=n!\binom{k}{m}k(k-1)\cdots (k-m+1)\frac{(n-k)\cdots (n-2k+m+1)}{n(n-1)\cdots (n-k+1)}.
    \end{align}
    Since $k<n/2$ we have that $\binom{k}{m}<k^4<\frac{n^4}{16}$. Notice that
    \begin{align}\label{i2}
    \frac{\frac{k(k-1)\cdots (k-m+1)(n-k)(n-k-1)\cdots (n-2k+m+1)}{k!}}{\frac{n(n-1)\cdots (n-k+1)}{k!}}=\frac{\binom{n-k}{k-m}}{\binom{n}{k}}.
    \end{align}
    We want an upper bound on $P_{n,m}$ so notice that $\binom{n-k}{k-m}=\binom{n-k}{n-2k+m}$. Since $n-2k+m<8$ we have that $\binom{n-k}{k-m}<n^8$. Note that $\binom{n}{\lfloor n/2\rfloor}\ge \frac{2^n}{n+1}$ is largest among the binomial coefficients $\binom{n}{i}$. So $$\binom{n}{k}=\binom{n}{\lfloor n/2\rfloor}\frac{\lfloor n/2\rfloor}{n-k}>\frac{1}{2}\binom{n}{\lfloor n/2\rfloor}>\frac{2^n}{2(n+1)}.$$ Now from (\ref{pnkm}) and (\ref{i2}) we get that $$P_{n,m}\le n!\frac{n^4}{16}\frac{n^82(n+1)}{2^n}=\frac{n^{12}(n+1)}{8}\frac{n!}{2^n}.$$ Now the total number of transversals with at most four cells in $A\times B$ is $P_{n,0}+P_{n,1}+\cdots +P_{n,4}\le \frac{5n^{12}(n+1)}{8}\frac{n!}{2^n}$ so taking $p(x)=\frac{5x^{12}(x+1)}{8}$ gives the desired result.
\end{proof}

We can see that in Theorem \ref{thm3} we could not get a better bound using the same method by being more careful with the inequalities. We do however think that $f(n)$ is much lower. On the other end Erdős and Selfridge \cite{erdos} proved that in a game where winning sets have size $n$ there must be at least $2^{n-1}$ winning sets if it is a first player win. Thus $f(n)\ge 2^{n-1}$. 
However, the example given by Erd\H{o}s and Selfridge
\cite{erdos} cannot be embedded into a family of
transversals. Indeed, their example is equivalent to the family of all sets consisting of a fixed
element $v$ and exactly one element from each of $n-1$ disjoint pairs. There are $2^{n-1}$
such winning sets and player 1 has a winning strategy. However, two distinct winning sets may have intersection of size $n-1$, which is impossible for two distinct transversals. The following result however shows that for the Maker-Breaker version of the transversal game exponentially many winning sets are indeed sufficient to ensure a Maker win.

\begin{thm}
    Let $n\ge 5$ be a positive integer. One may choose a family $\mathcal{F}$ of at most $10^4\cdot3^n$ transversals in an $n\times n$ board so that the Maker-Breaker transversal $\mathcal{F}$-game is a Maker win.
\end{thm}
\begin{proof}
  We will first prove that Maker can win even by playing second on all sufficiently large boards and then we will stack many such boards in our $n\times n$ board allowing Maker to play several games as the second player. First we show the following claim.

  \textbf{Claim.} For any $m\ge 5$ Maker can win the transversal $m$-game by playing second (where Breaker plays first).

  \textit{Proof of Claim.} We first show that if Maker can win on an $m\times m$ board then he can also win on an $(m+1)\times (m+1)$ board. Indeed this is a simple induction step. By permuting rows or columns we may assume that Breaker plays his first move in cell $(1,2)$. Now Maker plays in $(1,1)$ and he continues to play on the $m\times m$ sub-board consisting of the cells in the set $\{2,3,\ldots ,m+1\}^2$. He can make a transversal there giving him a transversal on the full board.

  To prove the claim we just need to show that Maker can win for $m=5$. We may assume that Breaker plays in $(1,1)$. Now Maker plays in $(1,2)$. If Breaker plays his second move anywhere in the row 1 or column 2 then Maker can win by forming a transversal in the board obtained by removing row 1 and column 2 where now he plays first. This gives Maker a transversal of the full board. If however Breaker does not play his second move in row 1 or column 2 then there are essentially two cases. By permuting some of the rows $2,3,4,5$ and or columns $3,4,5$ we may assume that Breaker's second move is on $(2,1)$ or $(2,3)$. We deal with these separately.

  Case 1: Breaker's second move is $(2,1)$. Now Maker plays in $(3,1)$. Now consider the sub-board $T$ with row 3 and column 1 removed. Maker has already played in the top-left of that board while Breaker has not played in those cells. By using the same strategy as in Theorem \ref{maintransversal} we know that Maker can win on a $4\times 4$ board by playing first so he can guarantee a transversal in $T$. Together with $(3,1)$ he can guarantee a transversal of the full board.

  Case 2: Breaker's second move is $(2,3)$. Now Maker plays in $(2,1)$. Similar to the previous case consider the board $T'$ with row 2 and column 1 removed. Maker has played in the top-left of $T'$ while Breaker has not played in $T'$. Now since Maker can win by playing first on a $4\times 4$ board, Maker can form a transversal on $T'$. Together with $(2,1)$ this gives a transversal of the full board. This completes the proof of the claim.
  \QED

  Now we return to the proof of the theorem. Let $k,r$ be integers such that $n=5k+r$, where $0\le r\le 4.$ Consider the sets $S_i=\{5i-4,5i-3,5i-2,5i-1,5i\}^2$ of cells on the $n\times n$ board for $1\le i\le k-1$. Let $S=\{5k-4,5k-3,\ldots , n\}^2$ and $V=(\cup_{i=1}^{k-1}S_i)\cup S$. Now let $\mathcal{F}$ be the set of all transversals only consisting of cells in $V$. We will show that Maker can win the transversal $\mathcal{F}$-game. Indeed Maker may play his first move arbitrarily and then play as the second player on each sub-board $S_1,S_2,\ldots S_{k-1},S$ separately. Wherever Breaker plays Maker will play on the same sub-board in his next move (if possible) using the winning strategy from the above claim. If Breaker does not play on one of the sub-boards Maker can just pretend that he did and proceed accordingly (extra moves clearly do not hurt Maker). Thus he can form transversals in each of the sub-boards $S,S_1,S_2,\ldots S_{k-1}$ giving a transversal in $\mathcal{F}$. This means that Maker can guarantee winning the transversal $\mathcal{F}$-game.

  All that remains is to count the number of transversals entirely contained in $V$. Every such transversal is composed of one transversal from each $S_i$ and $S$. This yields at most $(5!)^k\cdot 6\cdot 7\cdot 8\cdot 9<10^4\cdot 3^n$ transversals. This completes the proof.
\end{proof}

The relation between fast strategies for weak games and explicit winning strategies for strong games was studied in \cite{ferberhefetz} for perfect matchings on complete graphs. In our setting the board is bipartite and the winning family $\mathcal{F}$ may contain only some of its perfect matchings. Even though a Maker win in the Maker-Breaker version does not guarantee a player 1 win in the original `Maker-Maker' version, since the opponent may complete a transversal first, it seems that it should not be that much harder (possibly with the addition of some extra transversals) for player 1 to win. This leads us to the following conjecture.
\begin{conj}
    There is a constant $c>1$ such that $f(n)<c^n$ for large enough $n$.
\end{conj}

Now we move on to $g(n)$. For each positive integer $n>5$ define $m_n$ as follows. Let $2=k_0< k_1<k_2<k_3<k_4=n$ be integers. Define $p_i=\prod\limits_{j=k_i}^{k_{i+1}-1}j$ for $0\le i\le 3$. Now let $m_n$ be the maximum possible value of $\min\limits_{0\le i\le 3} p_i$ over all possible $k_1,k_2,k_3$. We will show that $g(n)\le n!-m_n+1$. It is not hard to see that $$m_n-1=(n!)^{1/4-O(1/n)}.$$
Indeed since each factor in the product $n!=2\cdot 3\cdot \ldots \cdot n$ is at most $n$ one may pick (for $n\ge 8$ ensuring $n!>n^4$) $k_1,k_2,k_3$ in that order such that $(n!)^{1/4}/n\le p_i\le (n!)^{1/4}$ for $0\le i\le 2$. 
Since $p_0p_1p_2p_3=(n-1)!$, we also have
\[
p_3=\frac{(n-1)!}{p_0p_1p_2}
\geq \frac{(n!)^{1/4}}{n}.
\]
Thus $m_n\geq (n!)^{1/4}/n$.
Since $m_n\geq 2$ for $n\geq 8$, we obtain
\[
m_n-1\ge \frac{(n!)^{1/4}}{2n}
\]
This will be sufficient since $(n!)^2=\prod_{i=1}^ni(n+1-i)\ge n^n$ and we also clearly have that $m_n\le (n!)^{1/4}$. Thus $m_n-1=(n!)^{1/4-O(1/n)}$,
as claimed.

\begin{thm}
    For all $n>5$ we have that $g(n)\le n!-m_n+1$.
\end{thm}
\begin{proof}
Suppose that $\mathcal{F}$ is a family of at least $n!-m_n+1$ transversals. Let $\mathcal{G}=\mathcal{F}^c$ be the family of the $m_n-1$ transversals not included in $\mathcal{F}$. Let $2=k_0<k_1<k_2<k_3<k_4=n$ be such that $\min\limits_{0\le i\le 3}p_i=m_n$ where $p_i=\prod\limits_{j=k_i}^{k_{i+1}-1}j$. We follow the same strategy as in the proof of Theorem \ref{maintransversal} except that we now show: after the first $n-1$ $X$ moves when $X$ occupies all but one cell of the diagonal, for each transversal $T\in \mathcal{G}$ at least four of the $n-1$ cells occupied by $X$ do not belong to $T$. As before when permuting rows and columns we change $\mathcal{G}$ accordingly. Now consider what happens in the inductive proof for $(*)$ in Theorem \ref{maintransversal} when $X$ has made $k+1$ moves and $O$ has made $k$ moves. If $O$'s next move is on $(a,b)$ where $a\le k+1$ or $b\le k+1$ then $X$ can move anywhere in $(k+2,c)$ where $c\ge k+3$. Then we can swap columns $k+2$ and $c$ and then swap column $c$ and $k+3$ to ensure $(*)$. So $X$ had $n-k-2$ choices in the same row. Similarly if $a,b>k+1$ if say $b>k+2$ then $X$ can play in $(a,c)$ for any $c\ge k+2$ with $c\neq b$. Then we can swap rows $a$ and $k+2$ and swap some columns to ensure $(*)$. Likewise if $b=k+2$ then we can swap columns $k+2$ and $k+3$ and do the same as when $b>k+2$. Either way there are always at least $n-k-2$ choices in the same row for $X$. With a similar analysis $X$'s second move has at least $n-2$ choices in the same row and certainly his first move has at least $n-1$ choices in the same row. Now let $\mathcal{G}_n=\mathcal{G}$. We will form $\mathcal{G}_{n-1},\mathcal{G}_{n-2},\ldots\mathcal{G}_{k_3}$ as follows. Inductively at move $k\le n-k_3$ $X$ will move in the cell $C$ which is one of the possible choices described above that is in the least number of transversals in $\mathcal{G}_{n+1-k}$. Then we let $\mathcal{G}_{n-k}$ be the set of all transversals in $\mathcal{G}_{n+1-k}$ containing $C$. Since there are at least $n-k$ choices in the same row we must have that $|\mathcal{G}_{n-k}|\le \frac{|\mathcal{G}_{n+1-k}|}{n-k}$. Thus $$|\mathcal{G}_{k_3}|\le \frac{m_n-1}{\prod\limits_{j=k_3}^{n-1} j}<1$$ so $\mathcal{G}_{k_3}$ is empty and hence no transversal in $\mathcal{G}$ will contain all cells of the first $n-k_3$ moves. Similarly we can ensure that for each $1\le i\le 3$ no transversal of $\mathcal{G}$ contains all cells where $X$ made all his moves from move $n-k_i+1$ to $n-k_{i-1}$ inclusive. Indeed for each interval we start from $\mathcal{G}$ again and reduce it as much as possible at every step. This completes $n-2$ $X$ moves and indeed $X$ can make his next move in accordance with Theorem \ref{maintransversal}. Thus after $n-1$ $X$-moves we can ensure that $|S\setminus T|\geq 4$ for each transversal $T\in\mathcal{G}$, where $S$ is the set of cells of the first $n-1$ $X$-moves. Since $|S|=n-1$ and $|T|=n$, we have
\[
|T\setminus S|=|S\setminus T|+1\geq 5.
\]
 Now $X$ follows the same strategy as in Theorem \ref{maintransversal}. Note that all threats that $X$ uses have at most four cells that are not part of the first $n-1$ $X$ moves so $X$ can win with the same strategy and guarantee that his winning transversal is in $\mathcal{F}$. This shows that indeed $g(n)\le n!-m_n+1$.
\end{proof}

Note that $g(n)>(n-2)!$ since the family consisting of $(n-2)!$ transversals containing both $(1,1)$ and $(2,2)$ is clearly not winning for $X$. Thus $g(n)$ is much bigger than $f(n)$. So how big is $g(n)$? It seems that we need to remove a large portion of transversals to stop $X$ from winning. We finish with the following conjecture.
\begin{conj}
    As $n\rightarrow \infty$ we have that $g(n)=o(n!)$.
\end{conj}
\section{Conflict of Interest}
The author declares no conflict of interest.

\bibliographystyle{plain}
\bibliography{main1}

\end{document}